%% file: main.tex
\pdfoutput=1

\input{./files/commands.tex}

\documentclass[journal]{IEEEtran}

\usepackage[]{graphicx}
\usepackage{amsmath,amssymb,amsfonts,mathrsfs}
\usepackage{hyperref}
\usepackage{cleveref}
\usepackage[caption=false,font=footnotesize]{subfig}

\crefname{equation}{eq.}{eqs.}
\crefname{figure}{fig.}{figs.}
\newtheorem{remark}{Remark}[section]
\hypersetup{colorlinks, breaklinks, allcolors=[rgb]{0.54,0,0.54}}

\begin{document}

% -> Comment it so that only cited references show. 
% -> Include it to see all the references on the bib files
% \nocite{*}

\title{A Boundary Integral Method for 3D Nonuniform Dielectric Waveguide Problems via the\\ Windowed Green Function}

\author{
  Emmanuel~Garza,
  Constantine~Sideris,
  and~Oscar~P.~Bruno
  \thanks{CS gratefully acknowledges support by the National Science Foundation (1849965, 2047433) and the Air Force Office of Scientific Research (FA9550-20-1-0087). CS and EG gratefully acknowledge support by the National Science Foundation under Grant 2030859 to the Computing Research Association for the CIFellows Project. OB gratefully acknowledges support from NSF under contracts DMS-1714169 and DMS-2109831, by AFOSR under contract FA9550-21-1-0373 and by the NSSEFF Vannevar Bush Fellowship under contract N00014-16-1-2808.}
  \thanks{E. Garza and C. Sideris are with the Department of Electrical and Computer Engineering, University of Southern California, Los Angeles, CA 90089 USA.}
  \thanks{O.P. Bruno is with the Department of Computing and Mathematical Sciences, California Institute of Technology, Pasadena, CA 91125 USA.}
}

\maketitle

\begin{abstract}
  This paper proposes an efficient boundary-integral based ``windowed Green function'' methodology (WGF) for the numerical solution of three-dimensional electromagnetic problems containing dielectric waveguides. The approach, which generalizes a two-dimensional version of the method introduced recently, provides a highly effective solver for general electromagnetic problems containing waveguides. In particular, using an auxiliary integral representation, the proposed method is able to accurately model incident mode excitation. On the basis of a smooth window function, the integral operators along the infinite waveguide boundaries are smoothly truncated, resulting in errors that decay faster than any negative power of the window size. 
\end{abstract}


\IEEEpeerreviewmaketitle

%%%%%%%%%%%%%%%%%%%%%%%%%%%%%%%%%%%%%%%%%%%%%%%%%%%%%%%%%%%%%%%%%%%%%%%%%%%
%% Introduction
%%%%%%%%%%%%%%%%%%%%%%%%%%%%%%%%%%%%%%%%%%%%%%%%%%%%%%%%%%%%%%%%%%%%%%%%%%%

\section{Introduction}
\label{sec:intro}

In view of their ability to efficiently guide transfers of electromagnetic energy, dielectric (open) waveguides play central roles in many engineered electromagnetic systems, including antenna feeds, coaxial cable transmission lines, optical fibers and nanophotonic devices, among many others~\cite{Koenderink2015,knox1976dielectric,itoh1977application,doerr2008dielectric,Snyder1983}. As discussed in what follows, however, the computational simulation of electromagnetic waveguides in three-dimensions (3D) has posed a number of significant challenges---mostly concerning accuracy and computational cost. Seeking to tackle this difficulty, this paper introduces an efficient boundary integral equation (BIE) Green-function based methodology for this problem. This method incorporates, as a central enabling element, a novel ``windowed Green function'' strategy (WGF)---previously demonstrated in the context of layered media in two and three dimensions, and for waveguide structures in two dimensions~\cite{Bruno2016,Perez-Arancibia2017,bruno2017windowed,sideris2019ultrafast,Garza2020}---to handle the truncation of BIE integration domains for fully-vectorial problems involving 3D dielectric waveguides with arbitrary shapes and configurations.

Most of the existing approaches for waveguide simulation rely on volumetric discretizations and approximation of the differential form of Maxwell's equations~\cite{bierwirth1986finite}, incorporating an absorbing condition, such as a Perfectly Matched Layer (PML)~\cite{Berenger1994,taflove2005computational}, to truncate the computational domain and, in particular, all of the semi-infinite waveguides present in the structure. The approach thus requires use of large 3D computational domains, and, consequently, large numbers of spatial unknowns. Further, evaluation of propagation over such large computational domains, which may amount to tens or even hundreds of wavelengths in relevant applications, can lead to significant dispersion errors by accumulation, over many wavelengths, of the errors inherent in the derivative approximations used in Finite Difference (FD) and Finite Element Method (FEM) at each discretization point. The time-domain versions of these methods that are often employed additionally suffer from dispersion and accuracy loss in the time variable. As a result, these approaches may require extremely fine meshes and associated high computing costs to maintain accuracy.

On the basis of the electromagnetic Green function, on the other hand, boundary integral methods (BIE) can lead to significantly reduced numbers of unknowns, since they only require discretization of the interfaces between different materials. Additionally, BIE methods are virtually free of numerical dispersion, in view of the Green function's ability to analytically (and, thus, without dispersion) propagate oscillatory fields to arbitrary distances. And, although in their straightforward implementations they lead to computing costs that grow quadratically with the discretization sizes (at least if iterative solvers such as GMRES~\cite{saad1986gmres} are used), BIE methods can be accelerated by a variety of techniques~\cite{Greengard1987,Greengard1998,Gumerov2004,Bruno2001b,Bruno2001}. Unfortunately, the application of BIE methods has almost exclusively been restricted to systems with bounded scatterers, in view of the lack of a suitable condition for termination of the infinite computational domain arising from infinite boundary interfaces.

Several effective approaches have recently emerged for the BIE waveguide-truncation problem, however. These include the WGF method in the frequency~\cite{Bruno2016,Perez-Arancibia2017,bruno2017windowed,sideris2019ultrafast,Garza2020} and time \cite{Labarca2019} domains, the closely related integral PML method~\cite{lu2018perfectly}, and a physically-motivated ``Surface Conductive Absorber" approach (SCA)~\cite{Zhang2011,zhang2010boundary}. Other approaches based on {\em volumetric} integral methods using adiabatic absorbers have also been demonstrated~\cite{tambova,Tambova2018}. The WGF method relies on a reformulation of the waveguide BIE equations, which results as the Green function is smoothly windowed to a compact support via multiplication by a ``slow-rise'' window function in the integration variable. This procedure, which only requires multiplication of the Green function by an adequate smooth window function, effectively screens out infinite boundary portions without introducing nonphysical reflection points on the boundary of the truncation domain, and it thus yields super-algebraic convergence---that is, asymptotically, faster convergence than any negative power of the window size. The integral PML contribution~\cite{lu2018perfectly} operates on the basis of the same principle as the WGF method: in the integral PML case the windowing effect is achieved by complexifying the spatial variables starting at a certain distance along the waveguide, which results in decaying exponentials in the Green function and thus provide the desired reflection-screening effect. The SCA approach~\cite{Zhang2011}, finally, incorporating a surface conductivity on the material interfaces sufficiently far along the waveguides, demonstrates decay rates of the order of $1/x^p$ with the distance $x$ along the waveguide, with e.g. $4\leq p\leq 8$ for screening of waveguide modes and with $p=1$ for screening of reactive fields; see~\cite[Figs. 9,~10]{Zhang2011}.

In this contribution, we generalize the 2D WGF method for waveguides~\cite{bruno2017windowed} to the fully-vectorial 3D counterpart. We show that our approach can handle both illuminating beams and point sources, as well as direct mode illumination. This paper is organized in the following way. In \cref{sec:bie} we describe the integral representation of the fields and the associated integral equations over the unbounded waveguide boundaries, making a distinction between the illuminating fields, with~\cref{sec:type-i} describing illuminating beams and point sources, and \cref{sec:type-ii} deals with direct mode sourcing. \Cref{sec:win} describes the ideas underlying the windowing of integral operators, and it presents our implementations for the two types of incidence considered in this paper. In particular \cref{sec:modes} describes the method we propose to accurately incorporate incident bound mode excitation---which, does not give rise to  undesirable radiating components. Finally, \cref{sec:num_ex} presents a variety of numerical examples that demonstrate the effectiveness and flexibility of the WGF method for geometrically complex waveguide structures.

%%%%%%%%%%%%%%%%%%%%%%%%%%%%%%%%%%%%%%%%%%%%%%%%%%%%%%%%%%%%%%%%%%%%%%%%%%%
%% BIE formulation
%%%%%%%%%%%%%%%%%%%%%%%%%%%%%%%%%%%%%%%%%%%%%%%%%%%%%%%%%%%%%%%%%%%%%%%%%%%

\section{Boundary Integral Formulation for Three-Dimensional Dielectric Waveguides}
\label{sec:bie}

Our approach is based on M\"{u}ller's frequency-domain integral formulation~\cite{muller2013foundations,yla2005well} at a given temporal frequency $\omega$, which follows from use of Green representation theorems~\cite{Nedelec2001}.  For clarity, in our description we consider a waveguide structure consisting of a single core (interior) region $\Omega_i$ (containing a material of permittivity $\perm_i$, magnetic permeability $\mu_i$ and spatial wavenumber $k_i$), and a single cladding (exterior) region $\Omega_e$ (containing material of permittivity $\perm_e$, magnetic permeability $\mu_e$ and spatial wavenumber $k_e$), as depicted in Fig.~\ref{fig:domains}.  The \emph{unbounded} interface between $\Omega_i$ and $\Omega_e$ and the corresponding normal pointing from the former to the latter are denoted by $\bdry$ and $\nn$, respectively. The wavenumber is related to the frequency and material properties by the relation $k = \omega \sqrt{\perm \mu}$. The generalization of these methods to structures including an arbitrary number of waveguides and dielectric materials does not pose major difficulties.

In this paper we consider the following two different types of incident excitation, namely:
\begin{itemize}
  \item \textbf{Type I excitation} (Beams or point source incidence): In this case we express the \emph{total} electric field $\Efield$ as the sum of the incident and scattered fields:
        \begin{align}
          \Efield =
          \begin{cases}
            \Efield_i + \Efield^\inc_i & \text{in } \Omega_i, \\
            \Efield_e + \Efield^\inc_e & \text{in } \Omega_e.
          \end{cases}
        \end{align}
  \item \textbf{Type II excitation} (Bound-mode incident field): The \emph{total} electric field $\Efield$ in this case is denoted as follows:
        \begin{align}
          \Efield =
          \begin{cases}
            \Efield_i & \text{in } \Omega_i, \\
            \Efield_e & \text{in } \Omega_e.
          \end{cases}
        \end{align}
\end{itemize}
The corresponding expressions from the magnetic field result from the replacement $\Efield\leftrightarrow \Hfield$.

\begin{remark}
  For both type I and type II problems, the electromagnetic field can be expressed in terms of Green function-based integral representations. The main difference in our treatment of these two cases is that for type II problems we use an integral representation of the \emph{total} field (including the incident bound mode), whereas for type I problems we use an integral representation of the \emph{scattered} field, such as is used for problems of scattering by bounded obstacles.
\end{remark}

In order to introduce the aforementioned integral representations, we consider the following vector potentials acting on a surface tangential density $\densa (\ney)$:
\begin{subequations}
  \label{eq:em-pot}
  \begin{align}
    \apot_\ell[\densa](\nex) \equiv & \; \nabla \times \int_\bdry
    G_\ell(\nex,\ney) \densa(\ney) \dsurf,                                                                     \\
    \bpot_\ell[\densa](\nex) \equiv & \; \nabla \times \nabla \times
    \int_\bdry G_\ell(\nex,\ney) \densa(\ney) \dsurf, \nonumber                                                \\
    =                               & \; k_\ell^2 \int_\bdry G_\ell(\nex,\ney) \densa(\ney) \dsurf + \nonumber \\
                                    & \; \; \nabla \int_\bdry G_\ell(\nex,\ney) \divsurf \densa(\ney) \dsurf.
  \end{align}
\end{subequations}
Here $G_\ell(\nex,\ney) = \exp{(ik_\ell|\nex-\ney|)} /(4\pi|\nex-\ney|)$ denotes the free-space Helmholtz Green function with wavenumber $k_\ell$. Subscript values $\ell = i$ and  $\ell =e$ indicate quantities corresponding to the interior and exterior domains, respectively. 

The tangential values of the field (3) for $\nex\in \bdry$  are given by the boundary integral operators
\begin{subequations}
  \label{eq:em-oper}
  \begin{align}
     & \soper_\ell [\densa] (\nex) \equiv - \nn(\nex) \times
    \int_\bdry G_\ell(\nex,\ney) \densa(\ney) \dsurf,              \\
     & \roper_\ell [\densa] (\nex) \equiv - \nn(\nex) \times \curl
    \int_\bdry G_\ell(\nex,\ney) \densa(\ney) \dsurf,              \\
     & \label{eq:toper}
    \toper_\ell [\densa] (\nex) \equiv - \nn(\nex) \times \nabla
    \int_\bdry G_\ell(\nex,\ney) \, \divsurf \densa(\ney) \dsurf,
  \end{align}
\end{subequations}
which are defined for $\nex\in\bdry$. For conciseness, in what follows we also use the weakly-singular operators
\begin{subequations}
  \label{eq:operdel}
  \begin{align}
    \roper^\Delta_\alpha & \equiv \frac{2}{\alpha_e+\alpha_i}
    \left( \alpha_e \roper_e -  \alpha_i \roper_i \right),             \\
    \koper^\Delta_\alpha & \equiv \frac{2i}{\omega(\alpha_e+\alpha_i)}
    \left[ (\toper_e - \toper_i) + (k_e^2 \soper_e -  k_i^2 \soper_i) \right],
  \end{align}
\end{subequations}
where the subindex $\alpha$ stands either for the dielectric constant symbol, $\alpha = \perm$, or the magnetic permeability symbol, $\alpha = \mu$. Finally, we call
\begin{subequations}
  \label{eq:dens}
  \begin{align}
    \densm_{i} \equiv - \nn \times \Efield_i \; \; & \text{and} \; \;
    \densj_{i} \equiv - \nn \times \Hfield_i \; \; \text{on} \; \bdry, \\
    \densm_{e} \equiv + \nn \times \Efield_e \; \; & \text{and} \; \;
    \densj_{e} \equiv + \nn \times \Hfield_e \; \; \text{on} \; \bdry.
  \end{align}
\end{subequations}
the tangential magnetic and electric surface currents $\densm_\ell$ and $\densj_\ell$ on the interior ($\ell = i$) and exterior ($\ell = e$) sides of $\bdry$.

Using these notations, the direct integral representations~\cite[theorem~5.5.1]{Nedelec2001} (cf. the Stratton-Chu formulas~\cite[theorem~6.7]{colton2013inverse}) for the exterior and interior fields are given by
\begin{subequations}
  \label{eq:rep}
  \begin{align}
    \apot_\ell[\densm_\ell] (\nex) + \frac{i}{\omega \perm_\ell} \bpot_\ell [\densj_\ell] (\nex) & =
    \begin{cases}
      \Efield_\ell (\nex) & \nex \in \Omega_\ell,      \\
      0                   & \nex \not \in \Omega_\ell,
    \end{cases}                                                                       \\
    \apot_\ell[\densj_\ell] (\nex) - \frac{i}{\omega \mu_\ell} \bpot_\ell [\densm_\ell] (\nex)   & =
    \begin{cases}
      \Hfield_\ell (\nex) & \nex \in \Omega_\ell,      \\
      0                   & \nex \not \in \Omega_\ell,
    \end{cases}
  \end{align}
\end{subequations}
with $\ell=e$ and $\ell = i$, respectively. The M\"{u}ller system of boundary integral equations results as  a linear combination of the tangential limiting values of the fields in~\cref{eq:rep} as $\nex \to \bdry$ from the exterior and interior domains~\cite{yla2005well,muller2013foundations,Hu2021}. The linear combination is selected in such a way that all the resulting integral operators are weakly-singular. In the following sections (\cref{sec:type-i,sec:type-ii}) we present the resulting system of integral equations for the two different types of source excitation.

\subsection{Type I Sources: Beams and Point Sources}
\label{sec:type-i}

In the case that the driving excitation equals either a planewave, a point source, or a beam (i.e. a superposition of planewaves or point sources), we obtain a system of integral equations that is entirely analogous to that obtained for a bounded dielectric obstacle---except, of course, that the integral operators are now defined over an \emph{unbounded} surface $\bdry$. The continuity of the tangent components of the electromagnetic fields tells us that the unknown density currents satisfy the relations
\begin{subequations}
  \label{eq:i-dens}
  \begin{align}
    \densm & \equiv \densm_i = - \densm_e + \nn \times \big( \Efield_i^\inc - \Efield_e^\inc \big), \\
    \densj & \equiv \densj_i = - \densj_e + \nn \times \big( \Hfield_i^\inc - \Hfield_e^\inc \big).
  \end{align}
\end{subequations}
The M\"{u}ller system  for the Type I unknowns $\densm$  and $\densj$ reads 
\begin{multline}
  \label{eq:ieq-sys}
  \begin{bmatrix}
    \idoper + \roper^\Delta_\perm & \koper^\Delta_\perm         \\
    -\koper^\Delta_\mu            & \idoper + \roper^\Delta_\mu
  \end{bmatrix}
  \begin{bmatrix}
    \densm \\
    \densj
  \end{bmatrix} = \\
  \begin{bmatrix}
    2 (\perm_e + \perm_i)^{-1} (\perm_e \Efield_e^\inc - \perm_i \Efield_i^\inc ) \times \nn \\
    2 (\mu_e + \mu_i)^{-1} (\mu_e \Hfield_e^\inc - \mu_i \Hfield_i^\inc)\times \nn
  \end{bmatrix}.
\end{multline}
It is easy to check that this system involves only weakly-singular kernels; see e.g.~\cite{yla2005well},~\cite[chapter~5]{Garza2020},~\cite{Hu2021}.

\subsection{Type II Sources: Bound Waveguide Modes}
\label{sec:type-ii}

Provided the refractive index in the core region is larger than that of the cladding, the structure can support bound modes and properly guide energy along the waveguide. As is known, an optical waveguide admits a finite number of guided bound modes of the form~\cite{Snyder1983}
\begin{align}
  \Efield_m(x,y,z) = \widetilde{\Efield}_m(x,y) \exp{(i k^m_z z)}
\end{align}
 (where it has been assumed the waveguide is aligned with the $z$ axis), with an analogous expression for $\Hfield$ (with the same propagation constant $k^m_z$). These modal fields decay exponentially fast away from the core and have a constant transverse profile.

Waveguide-mode excitation plays important roles in many realistic applications---since, e.g., it can be used to model fields incoming from other device components whose output is a waveguide mode. Consideration of waveguide-mode incident fields leads, in our formalism, to integral representations for the total fields, as indicated in what follows. These total-field integral representations, which incorporate integral representations for the incident bound mode as well as the radiation conditions for waveguide problems, take into account the directionality of incoming and outgoing modes~\cite{Nosich1994}.

Following~\cite{bruno2017windowed}, we assume the waveguide structure includes one or more ``semi-infinite waveguides'' (SIWs). A SIW is, simply, one half of an infinite waveguide (of arbitrary cross-section), as is obtained by cutting an ordinary infinite waveguide by a plane orthogonal to the optical axis. In addition to the SIWs, the structure may contain arbitrary dielectric structures as long as all such structures are contained within a bounded region: away from such a bounded region, only the SIWs break the homogeneity of space.

In view of these considerations, for the present incident-mode problems we define $\Omega^\inc$ as the union of all SIWs on which the incoming bound modes are given. Then, for both $\ell =i$ and $\ell = e$, the electric fields are characterized by the decomposition
\begin{align}
  \Efield_\ell =
  \begin{cases}
    \Efield_\ell^\scat + \Efield_\ell^\inc
                       & \text{in } \Omega_\ell \cap \Omega^\inc,                      \\
    \Efield_\ell^\scat & \text{in } \reals^3 \setminus (\Omega_\ell \cap \Omega^\inc),
  \end{cases}
\end{align}
with similar expressions for the magnetic fields $\Hfield_\ell$. Note that the scattered fields are discontinuous across some portions of the boundary of $\Omega^\inc$ which do not correspond to physical interfaces between different materials. The integral representation in~\cref{eq:rep} is still valid with the surface-current expressions in~\cref{eq:dens}, but, in the present incidence-mode case,  $\Efield_\ell$ ($\ell=i,e$) represents the total field, not just the scattered component of the field. We can then re-express the unknown surface currents in terms of the scattered fields, by utilizing the known incident densities:
\begin{subequations}
  \begin{align}
    \densm \equiv - \nn \times \Efield_i^\scat \; \;     & \text{and} \; \;
    \densj \equiv - \nn \times \Hfield_i^\scat \; \; \text{on} \; \bdry,    \\
    \densm^\inc \equiv - \nn \times \Efield_i^\inc \; \; & \text{and} \; \;
    \densj^\inc \equiv - \nn \times \Hfield_i^\inc \; \; \text{on} \; \bdry.
  \end{align}
\end{subequations}
Just as in the case of type I sources, we must enforce the continuity of the normal components of the fields, which imply that $\densm_e = - (\densm + \densm^\scat)$ and $\densj_e = - (\densj + \densj^\scat)$. Using the representation formula in~\cref{eq:rep} and with the same derivation as for type I problems, we obtain the M\"{u}ller system for type II incidence:
\begin{align}
  \label{eq:ieq-mode}
  \begin{bmatrix}
    \idoper + \roper^\Delta_\perm & \koper^\Delta_\perm         \\
    -\koper^\Delta_\mu            & \idoper + \roper^\Delta_\mu
  \end{bmatrix}
  \begin{bmatrix}
    \densm \\
    \densj
  \end{bmatrix} =
  - \begin{bmatrix}
    \idoper + \roper^{\Delta}_\perm & \koper^\Delta_\perm         \\
    -\koper^\Delta_\mu              & \idoper + \roper^\Delta_\mu
  \end{bmatrix}
  \begin{bmatrix}
    \densm^\inc \\
    \densj^\inc
  \end{bmatrix}.
\end{align}
In order to produce a physical solution (e.g., to avoid the trivial solution $[\densm, \densj] = -[\densm^\inc, \densj^\inc]$), the scattered currents $\densm$ and $\densj$ must satisfy the \emph{waveguide radiation conditions}~\cite[eq.~(36)]{Nosich1994}. Briefly, these radiation conditions, which are a generalization of the Silver-M\"{u}ller condition~\cite{muller2013foundations}, state that the scattered solution must only consist of \emph{outgoing} bound modes along the waveguide, and radiating fields away from the core region.

%%%%%%%%%%%%%%%%%%%%%%%%%%%%%%%%%%%%%%%%%%%%%%%%%%%%%%%%%%%%%%%%%%%%%%%%%%%
%% WGF
%%%%%%%%%%%%%%%%%%%%%%%%%%%%%%%%%%%%%%%%%%%%%%%%%%%%%%%%%%%%%%%%%%%%%%%%%%%
\section{Windowing of Integral Operators}
\label{sec:win}

The numerical solution of~\cref{eq:ieq-sys} over the unbounded waveguide boundary $\bdry$ requires careful consideration, in view of the slow decay rate---$\bigo(1/r)$---of the integrands in the associated integral operators. To do this, following~\cite{bruno2017windowed}, we: (1)~localize the scattering problem to a bounded region encompassing the nonuniform parts of the waveguide structure, (2)~use of a \emph{window} function to evaluate the slowly decaying, oscillatory integrals in the operators involving scattered currents, and (3) for type II problems (bound mode excitation) use of the representation formula for the modes to evaluate the operators involving the incident currents from the bound mode.\looseness = -1

The integrands in the integral operators~\cref{eq:em-oper}, which are defined over the \emph{unbounded} boundary $\bdry$, decay only as $\bigo(1/r)$ as $r\to\infty$, but they are also oscillatory, which makes them integrable. A direct truncation of the boundary $\bdry$ leads to extremely poor convergence---only $\bigo(1/\sqrt{A})$, where $A$ denotes the SIW-truncation length~\cite{bruno2017windowed,Monro2007}. However, \emph{super-algebraic} convergence---$\bigo(A^{-p})$ for any positive integer $p$---can be regained on the basis of the same truncation size $A$, simply by  multiplying the relevant integrands by a compactly supported, smooth ``window'' function. For strict super-algebraic convergence, the window function $w_A(d)$ must satisfy certain smoothness and ``slow-rise'' condition~\cite[section~III-B]{bruno2017windowed}, namely, that the windowing function to smoothly decay from $1$ to $0$ over the region $d \in (\alpha A, A) $ for some $\alpha \in (0, 1)$. High-order convergence can also be achieved using other suitable window functions based on functions that don't have strict compact support, but tend to zero exponentially, such as the error function or the hyperbolic tangent. A suitable choice of window function that satisfies the conditions for super-algebraic convergence is given by the expression
\begin{align}
  \label{eq:win}
  w_A(d) \equiv
  \begin{cases}
    1, & s(d) < 0,        \\
    \exp{\left(-2 \frac{\exp{(-1/|s(d)|^2)}}{|1-s(d)|^2}  \right)},
       & 0 \le s(d) \le 1 \\
    0, & s(d) > 1,
  \end{cases},
\end{align}
where $ s(d) = (|d|-\alpha A)/(A-\alpha A)$; throughout this paper we use the value $\alpha = 0.5$.

\begin{figure}[t!]
  \centering
  \includegraphics[width=\columnwidth]{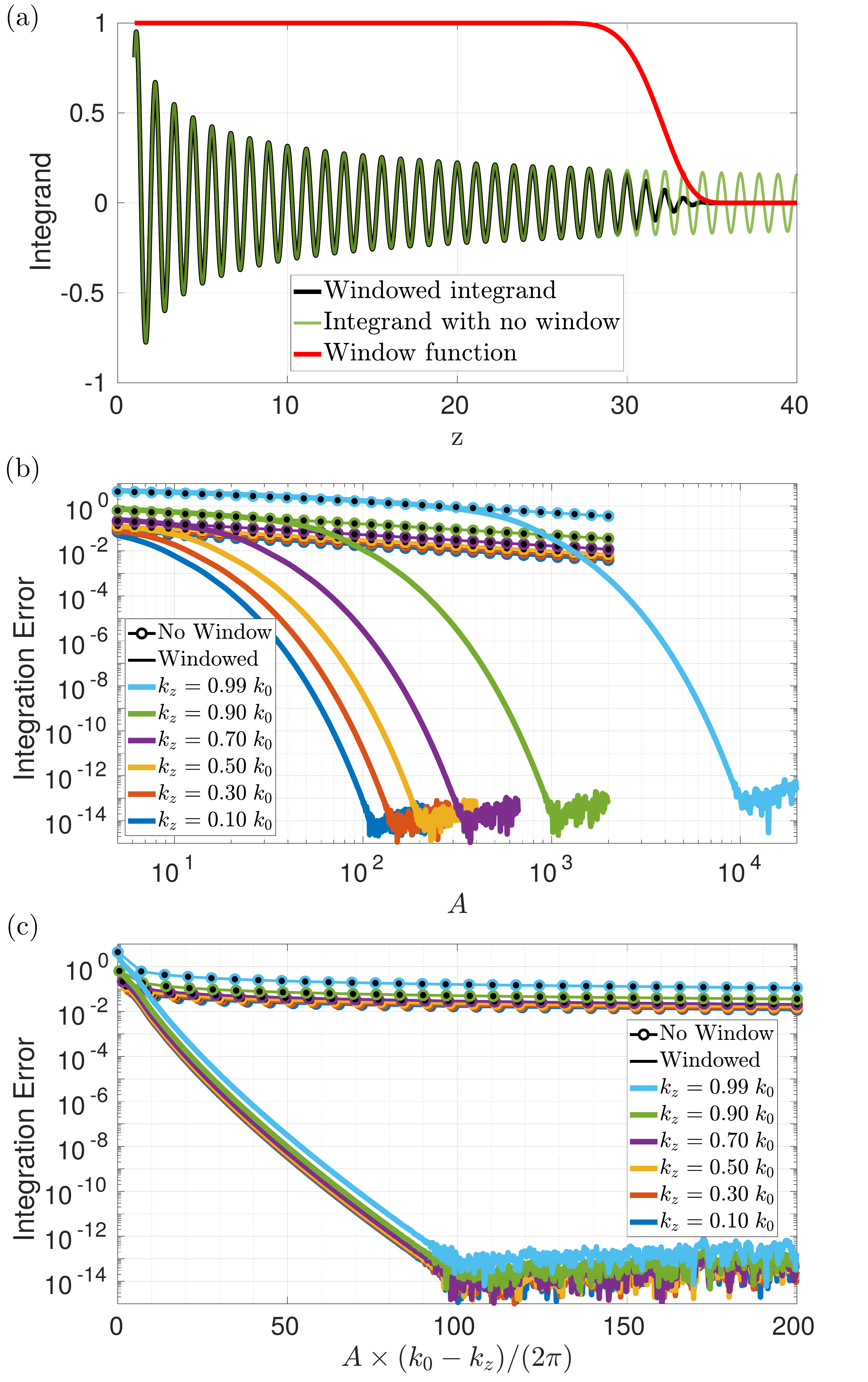}
  \caption{(a) Window function and windowed integrand for~\cref{eq:toy1,eq:toy2} for $k_z=0.1k_0$. (b) Log-log plot of convergence for~\cref{eq:toy1,eq:toy2} as we increase the window size $A$. (c) Semi-log plot of convergence as a function of the window size normalized by the \emph{net} wavelength.}
  \label{fig:toy}
\end{figure}

As a simple means to demonstrate the effectiveness of the windowing approach, for integrals similar to the ones needed to evaluate the operators in~\cref{eq:em-oper}, consider the one-dimensional integral given
\begin{align}
  \label{eq:toy}
  I = \int_1^\infty e^{i k_z z} \frac{e^{-i k_0 z}}{\sqrt{z}} \de z,
\end{align}
which is a good representation of the integrals needed for the integral operators in~\cref{eq:em-oper} due to the asymptotics of the Green's function, exemplified by the term $ e^{-k_0 z} / \sqrt{z}$, and the asymptotics of the current densities along the infinite boundaries, represented by the term $e^{i k_z z}$ in~\cref{eq:toy}. Note that for waveguide problems, these asymptotics are given by the modified radiation condition described in~\cite{Nosich1994}, which state that only bound modes propagate along the waveguide boundaries away from the nonuniform region. It can be shown, using integration by parts (see for example~\cite[section~2.2.3]{Monro2007},~\cite{Bruno2014}) that the approximations
\begin{align}
  \label{eq:toy1} I_{\text{tr}}(A) & = \int_1^A e^{i k_z z} \frac{e^{-i k_0 z}}{\sqrt{z}} \de z,        \\
  \label{eq:toy2} I_{w}(A)         & = \int_1^A w_A(z) e^{i k_z z} \frac{e^{-i k_0 z}}{\sqrt{z}} \de z,
\end{align}
have error rates (w.r.t. the value in~\cref{eq:toy}) given by
\begin{align}
  \bigo \left( \frac{1}{|k_0 -k_z|\sqrt{A}} \right)
  \; \; \; \text{and} \;\;\;
  \bigo \left( \frac{1}{|k_0-k_z|^p A^{p-\frac{1}{2}}} \right),
\end{align}
respectively, for any positive integer $p$.

\Cref{fig:toy} shows the convergence rates for the case $k_0 = 2\pi$ and several values of $k_z$ for both $I_\text{tr}(A)$ and $I_w(A)$. Two significant observations result from this simple yet illuminating example. On one hand, we see that using a window function $w_A$ can greatly improve the rate of convergence---so that certain integrals along infinite domains can be accurately truncated by just multiplying the integrand by $w_A$. On the other hand, we notice that, although the convergence rate is ``super-algebraic'' (faster than any power of $A$), if $k_z$ and $k_0$ are close to each other then impractically large values of $A$ may be necessary to achieve even modest accuracy. Hence, when the integrands might have vanishing oscillations---as it's the case of the right-hand-side of~\cref{eq:ieq-mode}---a special treatment is required.

In the context of the three-dimensional waveguide problems considered in this paper, we can construct a suitable window function, denoted by $\winoper(\nex)$ for $\nex \in \bdry$, on the basis of the 1D window function~\cref{eq:win}. Indeed, we first set $\winoper(\nex)$ to equal one when $\nex$ is not in any of the SIWs. Then, for $\nex$  in a SIW we set $\winoper$ equal to $w_A(d)$, where $d$ denotes the distance from $\nex$ to the plane perpendicular to the optical axis that determines the end of the SIW . To define this mathematically, let $M$ denote be the total number of SIWs, and let $\Omega_q^\text{SIW}$ denote the domain that encompasses the $q$-th SIW with it's origin located at $\oo_q$ and optical axis (pointing towards the infinite portion of the SIW) along the $\opax_q$ direction. Then, we have
\begin{align}
  \label{eq:win-op}
  \winoper(\nex) \equiv
  \begin{cases}
    w_A \left( \opax_q \cdot (\nex - \oo_q) \right),
       & \nex \in \bdry \, \cap \Omega_q^\text{SIW}, \\
    1, & \text{otherwise}.
  \end{cases}
\end{align}
This definition of the window function for our waveguide problems leads in a natural way to a truncated domain for which the window is non-zero, i.e. $\bdry^w = \bdry \, \cap \{\nex: \winoper(\nex) > 0 \}$. (For a different window function that does not have strict compact support, this definition of the truncated domain can be easily modified by taking the set for which $\winoper(\nex) > \eta$, for some small tolerance $\eta > 0$.)

%%%%%%%%%%%%%%%%%%%%%%%%%%%%%%%%%%%%%%%%%%%%%%%%%%%%%%%%%%%%%%%%%%%%%%%%%%%
%% Direct windowing
%%%%%%%%%%%%%%%%%%%%%%%%%%%%%%%%%%%%%%%%%%%%%%%%%%%%%%%%%%%%%%%%%%%%%%%%%%%
\subsection{Type I Incidence: Direct Windowing}
\label{sec:i-win}
For the sake of simplicity, we assume  incident fields for type I problems to have wavevectors with positive projections onto the optical axis $\opax_q$ of all the present SIWs. This condition guarantees that the net oscillation of the integrands in the operators of~\cref{eq:ieq-sys}---taking into account both the oscillation of the kernels and the solution currents---are bounded-away from zero, so that we can rely on windowing along the infinite boundary $\bdry$ to accurately truncate the simulation boundaries onto $\bdry^w$. In the case that this condition is not met, such vanishing oscillations can be treated by means of an approach similar to the one used in~\cite{Bruno2016,Perez-Arancibia2017} for layered media, or, simply, by re-scaling the infinite boundaries in terms of the net integrand wavenumber (net number of oscillations).

\begin{figure}[ht]
  \centering
  \includegraphics[width=1\columnwidth]{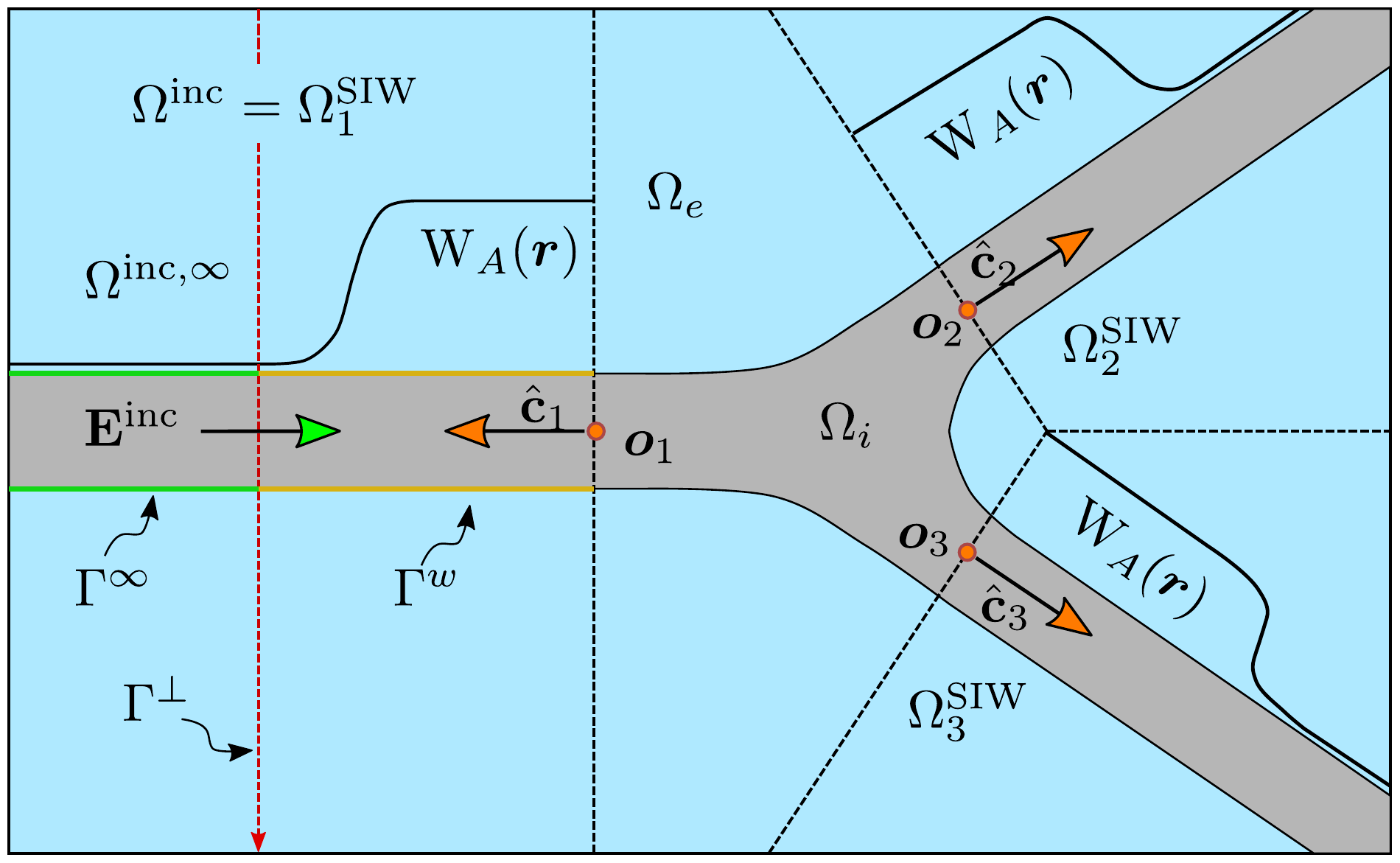}
  \caption{Diagram of the windowing strategy and auxiliary domains used for bound mode excitation.}
  \label{fig:domains}
\end{figure}

With the aforementioned assumption on the illuminating fields, we can directly window the integral operators on the left-hand side of~\cref{eq:ieq-sys} to obtain
\begin{multline}
  \label{eq:ieq-sys-w}
  \begin{bmatrix}
    \idoper + \roper^\Delta_\perm \winoper & \koper^\Delta_\perm\winoper           \\
    -\koper^\Delta_\mu \winoper            & \idoper + \roper^\Delta_\mu  \winoper
  \end{bmatrix}
  \begin{bmatrix}
    \densm \\
    \densj
  \end{bmatrix} = \\
  \begin{bmatrix}
    2 (\perm_e + \perm_i)^{-1} (\perm_e \Efield_e^\inc - \perm_i \Efield_i^\inc ) \times \nn \\
    2 (\mu_e + \mu_i)^{-1} (\mu_e \Hfield_e^\inc - \mu_i \Hfield_i^\inc)\times \nn
  \end{bmatrix}.
\end{multline}
This system of integral equations over the \emph{bounded} boundary $\bdry^w$ provides a super-algebraic approximation (w.r.t. the window size $A$) of the original, \emph{unbounded} problem in the region of interest, i.e. where the nonuniform part of the structure is present. The system in~\cref{eq:ieq-sys-w} can now be solved numerically using any bounded obstacle numerical method for integral equations.

\begin{remark} \label{rem:osc}
  In view of the definitions of the currents $\densj$ and $\densm$ in terms of the interior scattered fields in~\cref{eq:i-dens}, together with the waveguide radiation conditions discussed in~\cref{sec:type-ii}, the currents behave asymptotically as the tangential components of a superposition of outgoing bound modes. The $m$-th mode contribution along the $q$-th SIW contains a factor of $\exp{(+i k^m_q |\opax_q \cdot \ney|)}$ (where $k^m_q$ denotes the propagation constant of the $m$-th mode). Since the integral kernels oscillate with a factor of $\exp{(+ik|\nex-\ney|)}$, and since both $k^m_q$ and $k$ are positive, the product of the kernels and the currents result in non-vanishing oscillations as $\nex|_\bdry \to \infty$.
\end{remark}

%%%%%%%%%%%%%%%%%%%%%%%%%%%%%%%%%%%%%%%%%%%%%%%%%%%%%%%%%%%%%%%%%%%%%%%%%%%
%% Incident Modes
%%%%%%%%%%%%%%%%%%%%%%%%%%%%%%%%%%%%%%%%%%%%%%%%%%%%%%%%%%%%%%%%%%%%%%%%%%%

\subsection{Type II Incidence: Accurate Evaluation of Incident Modes}
\label{sec:modes}

In many instances, illuminating a waveguide with an incoming bound mode is desirable. Historically, the accurate sourcing of bound modes has been a non-straightforward matter, and alternative approximations are typically used---such as mode bootstrapping, illumination by Gaussian beams that approximate the mode, or by exciting the modes with point sources~\cite{taflove2005computational,Zhang2011}. These techniques, usually require either additional simulations, or large propagation distances for the incoming waves to shed away the undesired radiative or modal components, or the simulation is restricted to single-mode waveguides to avoid the excitation of spurious modes. However, it is highly advantageous to be able to directly source \emph{any} mode at will, thus avoiding unnecessary computation and potential errors.  With this goal in mind, this section proposes an integral equation methodology to accurately simulate the scattering of incident bound modes, so that any given incident bound mode for the relevant waveguides can be considered---on the basis of an auxiliary representation for the incident fields which, at minimal expense, incurs in errors that are exponentially small with regards to a certain truncation parameter related to the aforementioned auxiliary representation.

Windowing the integration domain in the same way as in~\cref{sec:i-win}, we can split the integrals on the right-hand side of~\cref{eq:ieq-mode} as the sum of two integrals---one over the bounded domain $\bdry^w$ and another over the unbounded domain $\bdry^\infty$ (so that $\bdry = \bdry^w \cup \bdry^\infty$). Let $\roper^{\Delta, w}_\alpha$ and  $\roper^{\Delta, \infty}_\alpha$ denote the operators in~\cref{eq:operdel} but with the integration domains in~\cref{eq:em-oper} substituted by $\bdry^w$ and $\bdry^\infty$, respectively. Using similar definitions  for the $\koper$ operator, we see that the right-hand side operators for the first equation in~\cref{eq:ieq-mode} are given by\looseness = -1
\begin{multline}
  \roper_\perm^{\Delta}[\densm^\inc](\nex) + \koper_\perm^{\Delta}[\densj^\inc](\nex) = (\roper_\perm^{\Delta,w} + \roper_\perm^{\Delta,\infty})[\densm^\inc](\nex) +  \\
  (\koper_\perm^{\Delta,w} + \koper_\perm^{\Delta,\infty})[\densj^\inc](\nex).
\end{multline}
In this expression, the terms involving integrals over $\bdry^w$ can be easily computed, since the integration domain is bounded. On the other hand, the integrands decay slowly over the infinite surface $\bdry^\infty$ for target points $\nex$ in the region of interest ($\nex \in \bdry^w$), and may have vanishing oscillations since the direction of the incoming mode is opposite to that of the oscillations from the integral kernels, much like the simple one-dimensional integral in~\cref{eq:toy}.

\begin{figure}[t!]
  \centering
  \includegraphics[width=0.96\columnwidth]{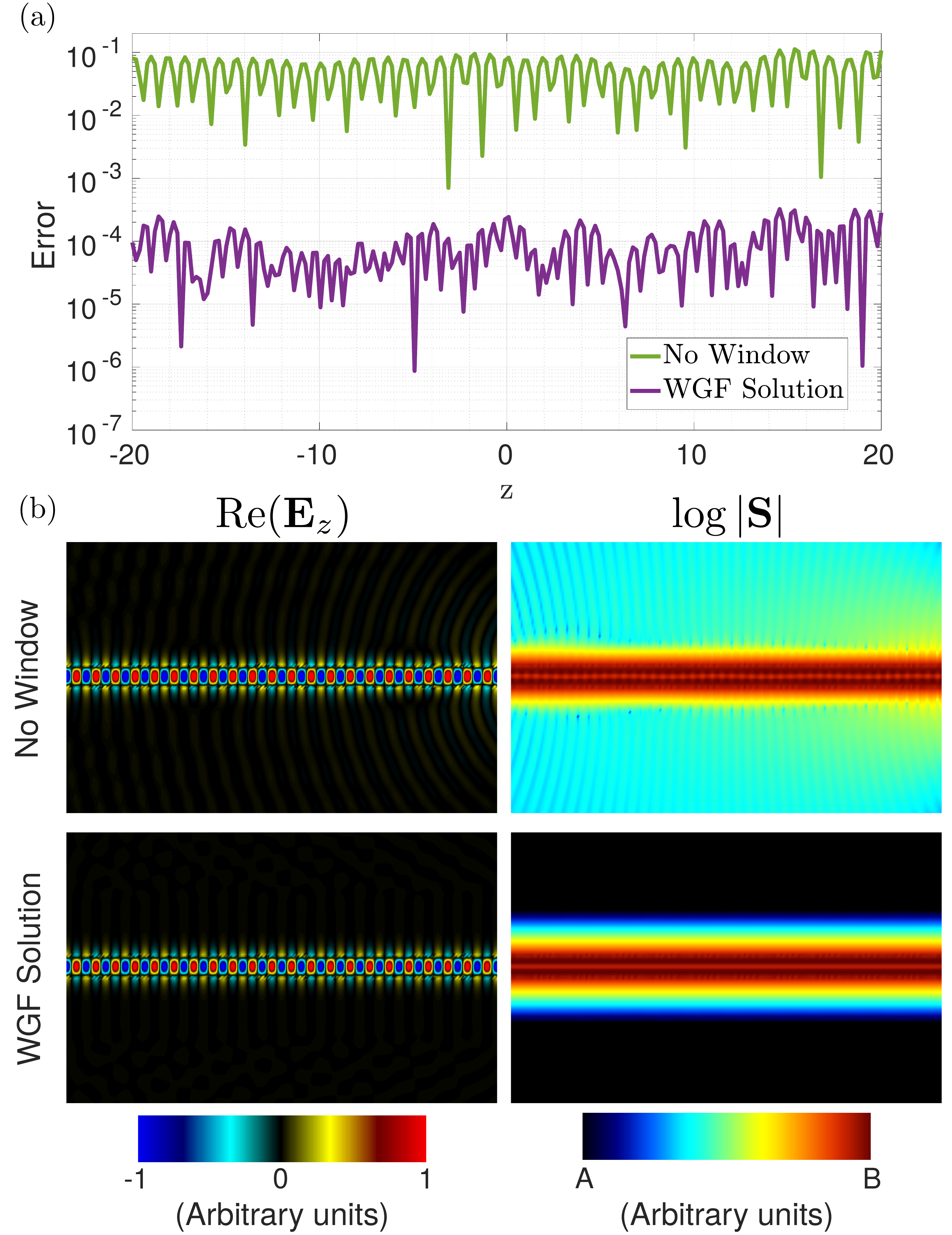}
  \caption{The error in $\Efield$ along the center of the waveguide is presented in (a). The top row of (b) shows the real part of the $z-$component of the electric field and the absolute value of the Poynting vector, in logarithmic scale, for the solution obtained without windowing the integral operators in~\cref{eq:ieq-mode-w}. The bottom row shows the same quantities, but with the WGF solution. The color scales are in arbitrary units and are the same color schemes that will be used in all figures. For this particular case, the color scale for $\log{|\mathbf{S}|}$ goes from $A=1\times10^{-6}$ to $B=10^{-1}$. \label{fig:exact} }
\end{figure}

To evaluate the integrals over $\bdry^\infty$ with minimal error, we introduce an auxiliary representation for the incident modes. Let $\bdry^\perp$ denote the orthogonal plane that cuts the incident SIW at the end of $\bdry^w$, and let $\Omega^{\inc,\infty}$ denote the portion of $\Omega^\inc$ that goes from $\bdry^\perp$ towards the infinite side of the SIW (see~\cref{fig:domains}). Using the representation theorems for the electromagnetic fields~\cite{Nedelec2001} we obtain
\begin{multline}
  \label{eq:rep-aux}
  (\apot^{\infty}_{\ell} + \apot^{\perp}_{\ell})[\densm^\inc_\ell] (\nex) +
  \frac{i}{\omega \perm_\ell} (\bpot^{\infty}_{\ell} + \bpot^{\perp}_{\ell})
  [\densj^\inc_\ell] (\nex) \\
  =
  \begin{cases}
    \Efield^\inc_\ell (\nex) & \nex \in \Omega^{\inc,\infty}_\ell,      \\
    0                        & \nex \not \in \Omega^{\inc,\infty}_\ell.
  \end{cases}
\end{multline}

For any point $\nex \in \bdry^w$ (which are outside of $\Omega^{\inc,\infty}$) we then have the relation
\begin{multline}
  \apot^{\infty}_{\ell}[\densm^\inc_\ell] (\nex) +
  \frac{i}{\omega \perm_\ell} \bpot^{\infty}_{\ell} [\densj^\inc_\ell] (\nex)
  =\\
  - \apot^{\perp}_{\ell}[\densm^\inc_\ell] (\nex) -
  \frac{i}{\omega \perm_\ell} \bpot^{\perp}_{\ell} [\densj^\inc_\ell] (\nex)
\end{multline}
which can be used to evaluate the right-hand side of~\cref{eq:ieq-mode}:
\begin{multline}
  \roper_\perm^{\Delta,\infty}[\densm^\inc](\nex) + \koper_\perm^{\Delta,\infty}[\densj^\inc](\nex) =  \\
  -\roper_\perm^{\Delta,\perp}[\densm^\inc](\nex) - \koper_\perm^{\Delta,\perp}[\densj^\inc](\nex).
\end{multline}
Using a representation for the field $\Hfield$  similar to the  one used in~\cref{eq:rep-aux} for the field $\Efield$ we can obtain the right-hand side~\cref{eq:ieq-mode} in terms of integrals along the bounded surface $\bdry^w$ and the unbounded boundary $\bdry^\perp$. The resulting WGF version of~\cref{eq:ieq-mode} reads
\begin{multline}
  \label{eq:ieq-mode-w}
  \begin{bmatrix}
    \idoper + \roper^\Delta_\perm \winoper & \koper^\Delta_\perm \winoper         \\
    -\koper^\Delta_\mu\winoper             & \idoper + \roper^\Delta_\mu \winoper
  \end{bmatrix}
  \begin{bmatrix}
    \densm \\
    \densj
  \end{bmatrix} = \\
  - \begin{bmatrix}
    \idoper + (\roper^{\Delta,w}_\perm - \roper^{\Delta,\perp}_\perm) & (\koper^{\Delta,w}_\perm-\koper^{\Delta,\perp}_\perm)         \\
    -(\koper^{\Delta,w}_\mu - \koper^{\Delta,\perp}_\mu)              & \idoper + (\roper^{\Delta,w}_\mu - \roper^{\Delta,\perp}_\mu)
  \end{bmatrix}
  \begin{bmatrix}
    \densm^\inc \\
    \densj^\inc
  \end{bmatrix}.
\end{multline}

\begin{remark} \label{rem:trick}
  The operators on the right-hand side of \cref{eq:ieq-mode-w} involve integrals that can be accurately computed. The integrals over the bounded surface $\bdry^w$ can be treated as in the bounded obstacle case. On the other hand, the integrands over $\bdry^\perp$, namely, the incident bound modes, decay \emph{exponentially} as the distance to the interface tends to infinity---along $\bdry^\perp$. This exponential decay allows us to evaluate integrals along $\bdry^\perp$ by simple truncation, while incurring only exponentially small truncation errors.
\end{remark}

%%%%%%%%%%%%%%%%%%%%%%%%%%%%%%%%%%%%%%%%%%%%%%%%%%%%%%%%%%%%%%%%%%%%%%%%%%%
%% Examples
%%%%%%%%%%%%%%%%%%%%%%%%%%%%%%%%%%%%%%%%%%%%%%%%%%%%%%%%%%%%%%%%%%%%%%%%%%%

\section{Numerical Examples}
\label{sec:num_ex}
\begin{figure}[t!]
  \centering
  \includegraphics[width=\columnwidth]{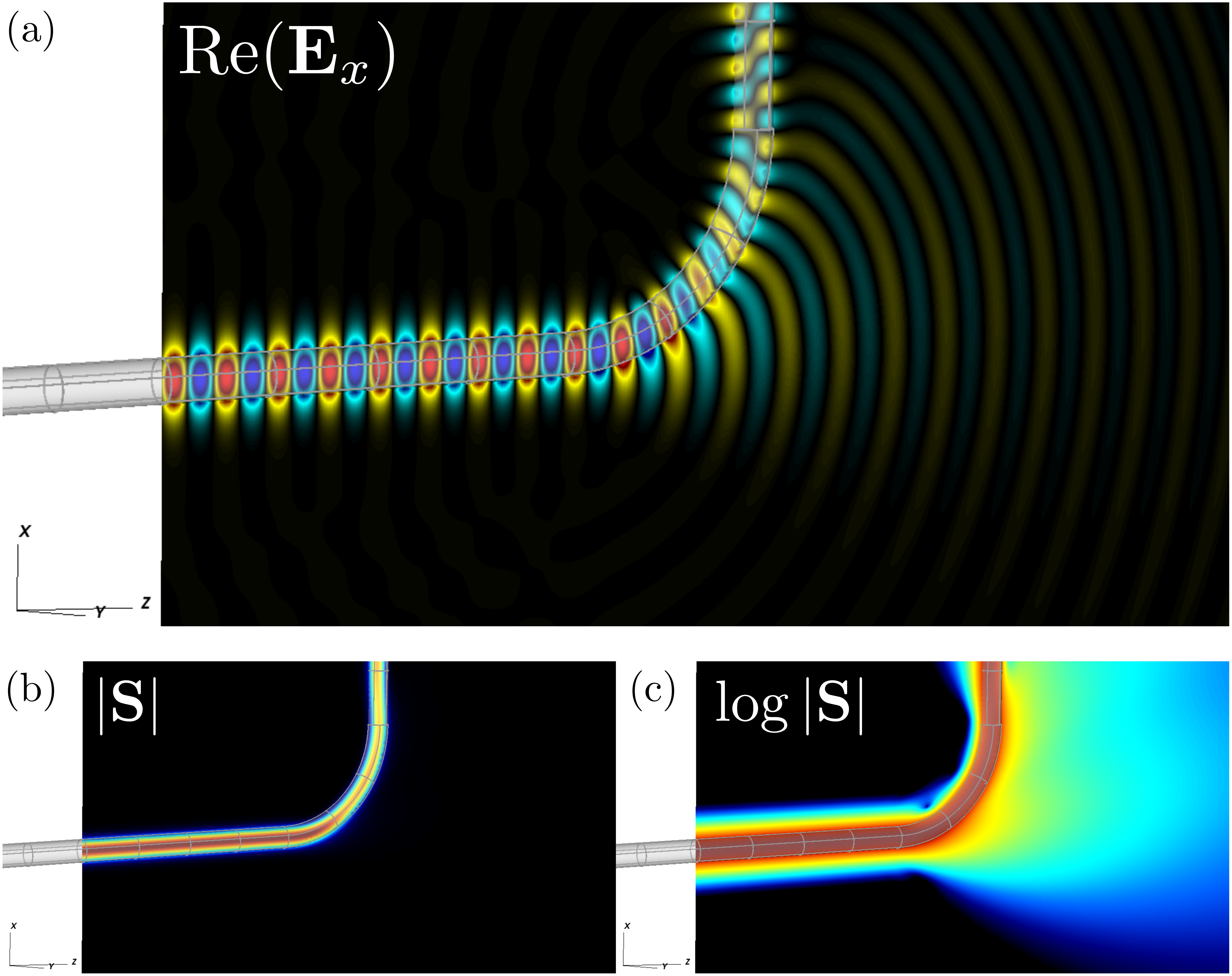}
  \caption{Mode propagation along a $90^\circ$ bend. The color scale for $|\mathbf{S}|$ and $\log{|\mathbf{S}|}$ goes from $A=5 \times10^{-6}$ to $B=3\times10^{-1}$.  \label{fig:bend}}
\end{figure}

Any numerical discretization for boundary integral methods can be used to discretize the WGF integral equations: in particular, Galerkin, Nystr\"{o}m or collocation method can be used. The illustrations presented in this paper were produced by means of the integration strategy presented in~\cite{bruno2018chebyshev,Garza2020,Hu2021,Garza2021}. This approach relies on discretization of boundaries by a set of non-overlapping quadrilateral curvilinear patches. The unknown current densities are discretized via Chebyshev polynomials on each surface patch. The far interactions are treated via Fej\'{e}r's quadrature, and the near-singular and singular interactions are precomputed using a high-order quadrature rule for the weakly-singular integrals. Then, the system of integral equations is solved iteratively via GMRES.

\begin{figure}[t!]
  \centering
  \includegraphics[width=0.8\columnwidth]{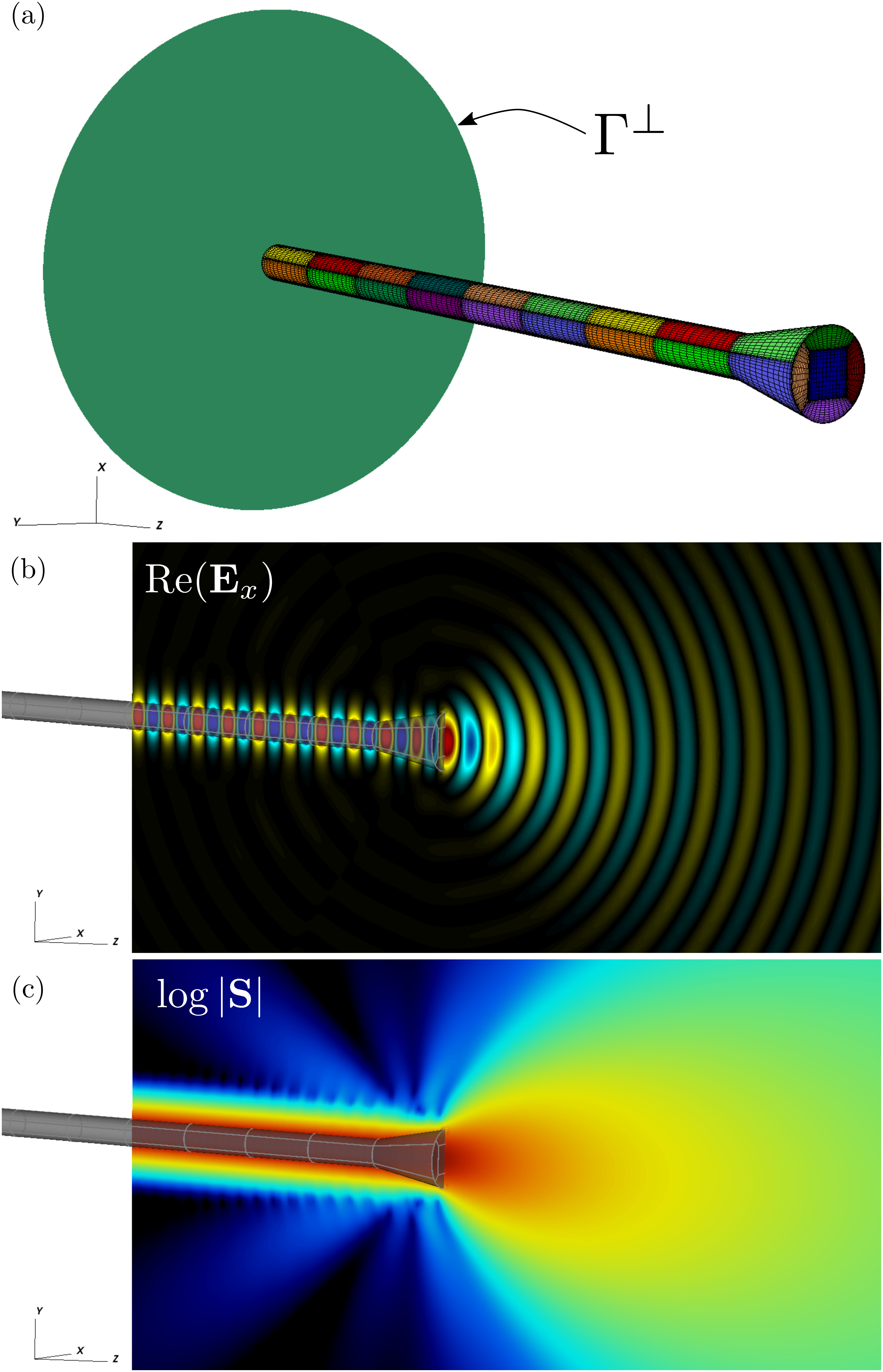}
  \caption{Modeling of a dielectric antenna by a terminated waveguide. (a) shows the discretization of the structure and the $\bdry^\perp$ used to evaluate the incident mode contribution. The color scale for (c) goes from $A=1 \times10^{-6}$ to $B=3\times10^{-1}$.  \label{fig:antenna}}
\end{figure}

This section presents several numerical examples that illustrate the accuracy and applicability of the proposed WGF approach. Our first example concerns  a uniform circular waveguide, which allows us to provide comparisons with the analytical solution of a mode propagating unperturbed along the waveguide. In our second example we then consider a circular waveguide with a 90$^\circ$ bend. In a third example, we consider a dielectric antenna---e.g. an open waveguide with a termination---demonstrating the applicability of the method for this class of devices. The last two examples, finally, concern illumination of waveguide structures by beams; the first one is that of an elliptical cross section waveguide, and the last one is that of a complex nanophotonic waveguide with multiple output waveguides. For all problems in this section we have assumed an exterior refractive index of $1.0$ and a interior refractive index of $1.47$ (SiO$_2$), with the exception of the last one for which the core is of silicon ($n_\text{Si}= 3.47$) and the cladding is of SiO$_2$.

\begin{figure}[t!]
  \centering
  \includegraphics[width=0.96\columnwidth]{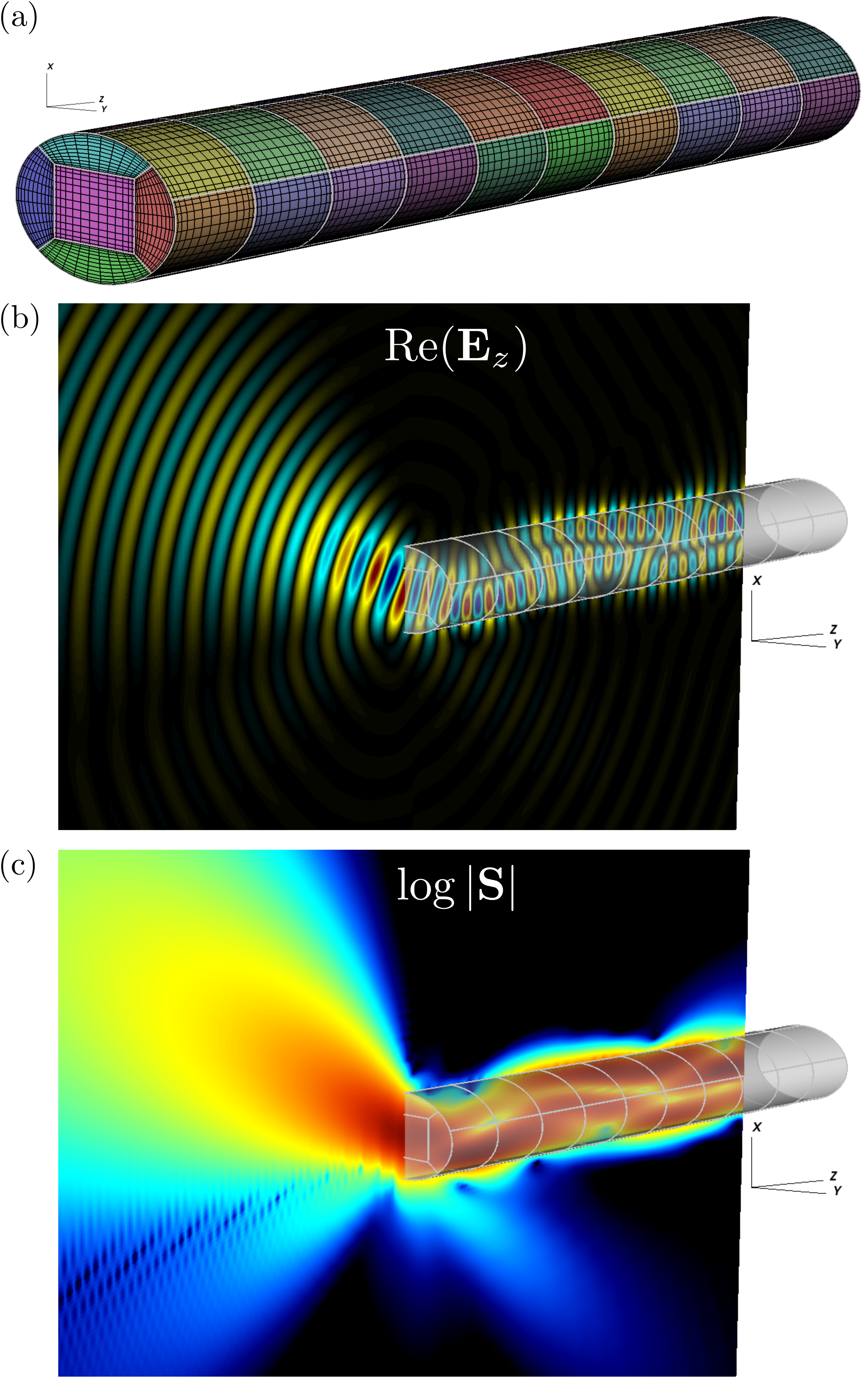}
  \caption{Illumination of an elliptical waveguide by an electromagnetic beam. In this mode launching problem, several modes get excited, and the simulation shows the ``bouncing'' of the trapped fields inside the waveguide. The color scale for (c) goes from $A=5 \times10^{-6}$ to $B=6 \times10^{-1}$. \label{fig:ell-ill}}
\end{figure}

\begin{figure}[ht]
  \centering
  \includegraphics[width=1\columnwidth]{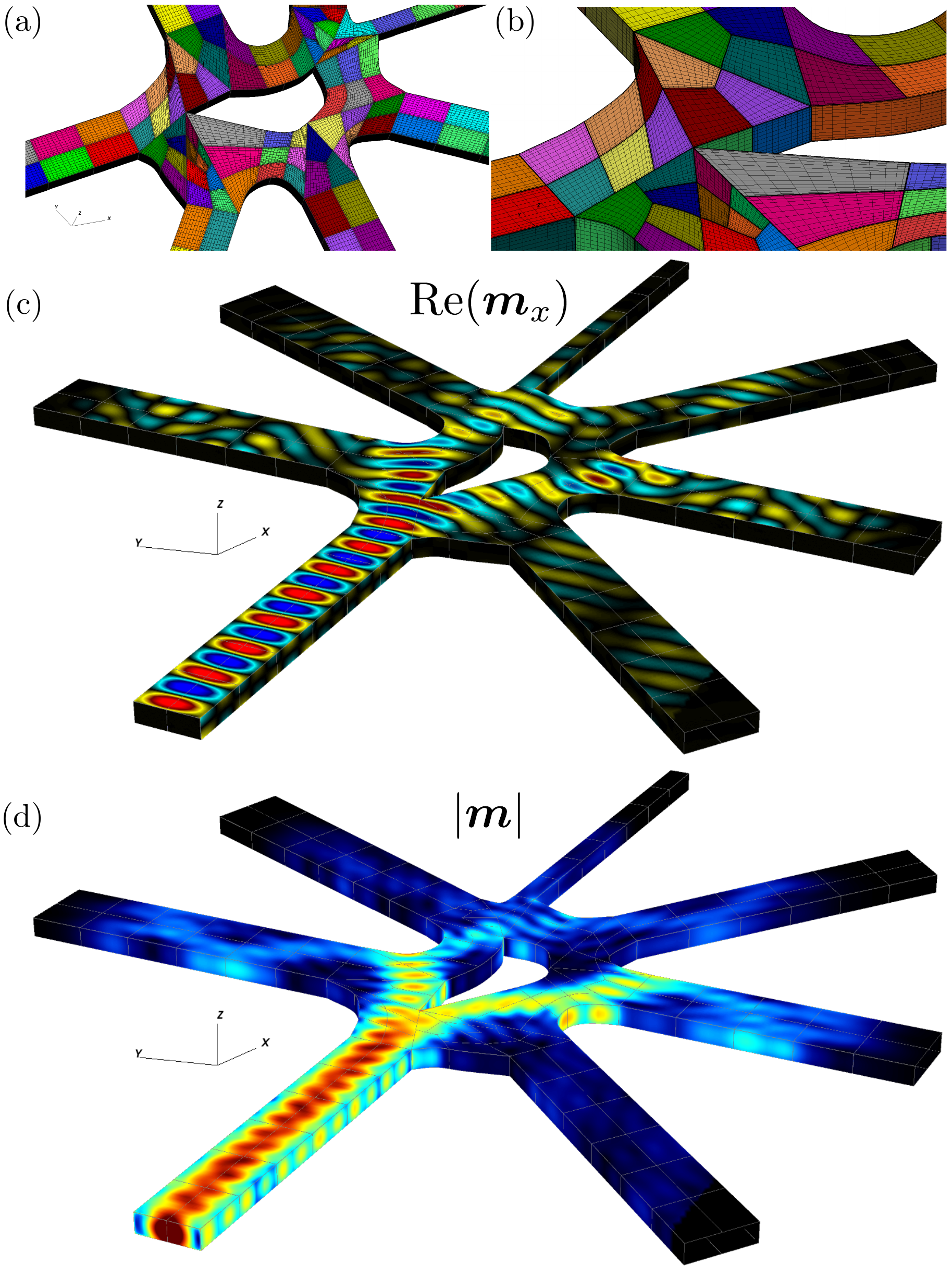}
  \caption{Scattering by a complex silicon waveguide structure. (a) and (b) show the transfinite patches used to discretize the geometry. Plots for the real part and the absolute value of the magnetic currents $\densm$ are shown in (c) and (d), respectively.\label{fig:multi}}
\end{figure}

%
% Uniform
%
In the case of a perfectly uniform waveguide with a circular cross section, the modal fields may be expressed in terms of Bessel functions~\cite{Snyder1983,Jackson1999}. \Cref{fig:exact} displays the corresponding numerical solution of~\cref{eq:ieq-mode-w}. In particular, \Cref{fig:exact}~(a) presents errors with respect to the exact mode solution for the modal fields along the line at the center of the waveguide. Errors obtained with and without use of a window function are presented in purple and green, respectively, on the basis of the same spatial discretization; the effect of the window function can easily be appreciated.  \Cref{fig:exact}~(b) displays the real part of the $z$-component of the electric field (left row) and the absolute value of the Poynting vector ($\mathbf{S} = \frac{1}{2} \Efield \times \Hfield^*$~\cite{Jackson1999}) in $\log$ scale, also for both, a solution obtained with no window and for the WGF solution. For this particular case, where a discretization of $18 \times 18$ ($\sim9$ points per mode's wavelength) is used, accuracies of order $10^{-4}$ are achieved by the WGF method. The un-windowed solution matches some qualitative features of the field, but it is otherwise clearly incorrect: the field values at the center of the waveguide are off (\cref{fig:exact}~(a)), and \cref{fig:exact}~(b) displays incorrectly curved wave fronts which arise as the physical (straight) wavefronts are superimposed to the nonphysical reflections arising from the truncation edge, as well as large mismatches in the energy flow.

%
% Bend
%
\Cref{fig:bend} presents results of a simulation of a circular waveguide with a $90^\circ$ bend illuminated by a bound mode. \Cref{fig:bend}~(a) plots the $x$-component of the electric field, and shows that the mode is mostly preserved across the curved region: the discrepancy in character of the field shown in the vertical and horizontal sections is only a matter of appearance, reflecting the fact that the $x$-component of the field is the tangential component in the vertical section, but it is the normal component in the horizontal section. In \cref{fig:bend}~(b) and (c) the absolute value of the Poynting vector is shown both in linear and logarithmic scales, showing that some of the energy radiates away from the waveguide as a result of the bend.

%
% Antenna
%
\Cref{fig:antenna} displays the fields resulting in a terminated open waveguide illuminated by a bound mode, i.e., a dielectric antenna. \Cref{fig:antenna}~(a) presents the discretization used as well as the auxiliary boundary $\bdry^\perp$ utilized in this case to evaluate the right-hand of~\cref{eq:ieq-mode-w}. \Cref{fig:antenna}~(b)~and~(c) present the real part of $\Efield_x$ and the absolute value of the Poynting vector (in logarithmic scale), displaying, in particular, the near-field pattern and the energy radiation away from the waveguide termination region.

%
% Elliptical
%
\Cref{fig:ell-ill} then present the fields obtained for a waveguide with elliptical cross-section illuminated by an electromagnetic beam~\cite{Chen2002}. In particular,  \cref{fig:ell-ill}~(b) and~(c) show that a significant fraction of the energy is coupled to the waveguide. The logarithmic scale in (c), in turn, reveals that parts of the beam also reflect backward,  and some of it passes through the structure. The problem considered in this example is representative of the ``mode launching'' problem, for which one illuminates a waveguide in order to couple energy into the waveguide modes. Note that although some pseudo-analytical methods exist for mode launching~\cite{Snyder1969a} they tend to rely on approximations, and, thus, the computational mode-launching capability could prove valuable. In particular, the WGF approach fully accounts, within the prescribed numerical error, for the transmitted and reflected fields.

%
% Photonic
%
Our last example, which is presented  in \cref{fig:multi}, concerns a complex silicon structure, where all of the SIWs have rectangular cross section. This example demonstrates the applicability of the WGF method to complex 3D nanophotonic structures. In order to handle the challenging geometrical features in this structure we utilized the CAD capabilities provided by gmsh~\cite{geuzaine2009gmsh} to create the necessary transfinite patches shown in \cref{fig:multi}~(a) and (b). The waveguide structure is illuminated by a Gaussian beam at the input port, and the solution to the magnetic current densities are shown in (c) and (d). \looseness = -1

%%%%%%%%%%%%%%%%%%%%%%%%%%%%%%%%%%%%%%%%%%%%%%%%%%%%%%%%%%%%%%%%%%%%%%%%%%%
%% Conclusions
%%%%%%%%%%%%%%%%%%%%%%%%%%%%%%%%%%%%%%%%%%%%%%%%%%%%%%%%%%%%%%%%%%%%%%%%%%%

\section{Conclusions}
We have presented a fully vectorial numerical method for the solution of complex 3D electromagnetic problems including waveguide structures on the basis of a windowed Green function boundary integral method. In this approach, the relevant integral operators---which are initially posed over the infinite boundaries of the waveguides---can be accurately evaluated, with super-algebraically small errors, by multiplying the integrands by a smooth window function and truncating the integration domain to the region where the window function does not vanish.  In particular, we showed that incident mode excitation can be accurately incorporated by means of an auxiliary representation, which transform challenging integrals along the infinite boundary of the waveguide carrying the incident mode, onto integrals with exponentially decaying integrands. The ideas presented here are independent of the numerical discretization of the integral operators, and can indeed be used, in particular, in conjunction with any Nystr\"{o}m, or standard method-of-moments approach.

%%%%%%%%%%%%%%%%%%%%%%%%%%%%%%%%%%%%%%%%%%%%%%%%%%%%%%%%%%%%%%%%%%%%%%%%%%%
%% References
%%%%%%%%%%%%%%%%%%%%%%%%%%%%%%%%%%%%%%%%%%%%%%%%%%%%%%%%%%%%%%%%%%%%%%%%%%%

\bibliographystyle{IEEEtran}
\bibliography{IEEEabrv,bibliography}

%%%%%%%%%%%%%%%%%%%%%%%%%%%%%%%%%%%%%%%%%%%%%%%%%%%%%%%%%%%%%%%%%%%%%%%%%%%
%% IEEE Biographies
%%%%%%%%%%%%%%%%%%%%%%%%%%%%%%%%%%%%%%%%%%%%%%%%%%%%%%%%%%%%%%%%%%%%%%%%%%%

% Insert where needed to balance the two columns on the last page with
% biographies
\newpage

\end{document}

%% file: files/commands.tex
\newcommand{\bdry}{\Gamma}
\newcommand{\curl}{\nabla \times}
\newcommand{\dsurf}{\text{d} \sigma(\ney)}
\newcommand{\divsurf}{\text{\textnormal{div}}_\bdry}

\newcommand{\apot}{\mathscr{A}}
\newcommand{\bpot}{\mathscr{B}}

\newcommand{\soper}{\operatorname{S}}

\newcommand{\toper}{\operatorname{T}}
\newcommand{\roper}{\operatorname{R}}

\newcommand{\koper}{\operatorname{K}}

\newcommand{\idoper}{\operatorname{I}}
\newcommand{\winoper}{\operatorname{W}_A}

\newcommand{\densj}{\boldsymbol{j}}
\newcommand{\densm}{\boldsymbol{m}}
\newcommand{\densa}{\boldsymbol{a}}

\newcommand{\inc}{\text{inc}}
\newcommand{\scat}{\text{scat}}
\newcommand{\Efield}{\mathbf{E}}
\newcommand{\Hfield}{\mathbf{H}}

\newcommand{\de}{\text{d}}
\newcommand{\bigo}{\mathcal{O}}
\newcommand{\reals}{\mathbb{R}}
\newcommand{\nn}{\boldsymbol{n}}
\newcommand{\nex}{\boldsymbol{r}}
\newcommand{\ney}{\boldsymbol{r}'}

\newcommand{\oo}{\boldsymbol{o}}
\newcommand{\opax}{\hat{\boldsymbol{c}}}

\newcommand{\perm}{\varepsilon}